\documentclass[12pt,thmsa]{article}
\usepackage{amsfonts}

\input{tcilatex}
\begin{document}

\author{Karl-Georg Schlesinger \qquad \\
Institute for Theoretical Physics\\
Technical University of Vienna\\
Wiedner Hauptstr. 8-10\\
A-1040 Wien, Austria\\
e-mail: kgschles@esi.ac.at}
\title{A duality Hopf algebra for holomorphic $N=1$ special geometries}
\date{}
\maketitle

\begin{abstract}
We find a self-dual noncommutative and noncocommutative Hopf algebra $%
\mathcal{H}_F$ acting as a universal symmetry on the modules over inner
Frobenius algebras of modular categories (as used in two dimensional
boundary conformal field theory) similar to the Grothendieck-Teichm\"{u}ller
group $GT$ as introduced by Drinfeld as a universal symmetry of
quasitriangular quasi-Hopf algebras. We discuss the relationship to a
similar self-dual noncommutative and noncocommutative Hopf algebra $\mathcal{%
H}_{GT}$, previously found as the universal symmetry of trialgebras and
three dimensional extended topological quantum field theories. As an
application of our result, we get a transitive action of a sub-Hopf algebra $%
\mathcal{H}_D$ of $\mathcal{H}_{GT}$ on the relative period matrices of
holomorphic $N=1$ special geometries, i.e. $\mathcal{H}_D$ appears as a kind
of duality Hopf algebra for holomorphic $N=1$ special geometries.
\end{abstract}

\section{Introduction}

In \cite{Dri} Drinfeld has introduced the so called
Grothendieck-Teichm\"{u}ller group $GT$ as the universal symmetry of
quasitriangular quasi-Hopf algebras. We have shown in two previous papers (%
\cite{Sch 2002a}, \cite{Sch 2002b}) that one can extend this approach to get
a self-dual, noncommutative, and noncocommutative Hopf algebra $\mathcal{H}%
_{GT}$ acting as the universal symmetry on three dimensional extended
topological quantum field theories in the sense of \cite{KL}. We have also
shown there that these symmetry considerations lead to a very restrictive
stability property of three dimensional extended topological quantum field
theories which we called ultrarigidity, there.

The three dimensional extended topological quantum field theories of \cite
{KL} are motivated by the wish to describe two dimensional boundary
conformal field theories in a three dimensional topological way (much the
same way one can get the two dimensional WZW model in the bulk from three
dimensional Chern-Simons theory). There is another purely algebraic
formulation which directly applies to the full case of two dimensional
boundary conformal field theories (see \cite{FS}, \cite{FRS 2001}, \cite{FRS
2003}). There, two dimensional boundary conformal field theories with given
chiral data are formulated as modular categories $\mathcal{C}$ together with
a module $M$ over an inner Frobenius algebra of $\mathcal{C}$. We show that
one can, again, introduce a self-dual, noncommutative, and noncocommutative
Hopf algebra $\mathcal{H}_F$ as the universal symmetry of such pairs $\left( 
\mathcal{C},M\right) $. We discuss the relationship between $\mathcal{H}_F$
and $\mathcal{H}_{GT}$.

Finally, in section 3, we will apply our result to derive a transitive
action of a sub-Hopf algebra $\mathcal{H}_D$ of $\mathcal{H}_{GT}$ on the
relative period matrices as they are used in \cite{LMW 2002a}, \cite{LMW
2002b} to describe holomorphic $N=1$ special geometries. In this sense, $%
\mathcal{H}_D$ can be seen as a duality Hopf algebra for holomorphic $N=1$
special geometries.

\bigskip

\section{The Hopf algebra $\mathcal{H}_F$}

Let $\mathcal{C}$ be a modular category, i.e. a $\Bbb{C}$-linear, rigid,
abelian, braided monoidal category (observe that we include the condition
that $\mathcal{C}$ be monoidal into the definition of a modular category,
different variants of the precise definition are used in the literature). An
object $A$ in $\mathcal{C}$ is called an algebra object in $\mathcal{C}$ if
there exists a morphism 
\[
m:A\otimes A\rightarrow A
\]
in $\mathcal{C}$ which satisfies associativity and if there exists another
morphism $\eta $ in $\mathcal{C}$ from the unit object in $\mathcal{C}$ to $A
$ which satisfies the commutative diagrams of a unit object with respect to $%
m$. Dually one can define the notion of a coalgebra object $\left( A,\Delta
,\varepsilon \right) $ in $\mathcal{C}$ with $\Delta $ a coassociative
coproduct and $\varepsilon $ a counit. One calls $\left( A,m,\eta ,\Delta
,\varepsilon \right) $ a Frobenius algebra in $\mathcal{C}$ if $\left(
A,m,\eta \right) $ and $\left( A,\Delta ,\varepsilon \right) $ give an
algebra object and a coalgebra object in $\mathcal{C}$, respectively, and we
have 
\[
\left( id_A\otimes m\right) \circ \left( \Delta \otimes id_A\right) =\Delta
\circ m=\left( m\otimes id_A\right) \circ \left( id_A\otimes \Delta \right) 
\]
A (left) module over an algebra object in $\mathcal{C}$ is an object $M$ in $%
\mathcal{C}$ together with a morphism 
\[
\varrho :A\otimes M\rightarrow M
\]
such that 
\[
\varrho \circ \left( m\otimes id_M\right) =\varrho \circ \left( id_M\otimes
\varrho \right) 
\]
and 
\[
\varrho \circ \left( \eta \otimes id_M\right) =id_M
\]
In \cite{FS}, \cite{FRS 2001}, and \cite{FRS 2003} modules over Frobenius
algebras in modular categories have been shown to give an algebraic
formulation for boundary conditions of two dimensional conformal field
theories with given chiral data (for this formulation only so called
symmetric special Frobenius algebras appear; we refer the reader to the
literature cited above for this notion).

On the other hand, in \cite{KL} three dimensional extended topological field
theories (see there for the definition) have been introduced motivated by
the wish to give a three dimensional topological formulation of at least a
class of two dimensional boundary conformal field theories. We have shown in 
\cite{Sch 2002b} that any such three dimensional topological quantum field
theory uniquely determines a so called trialgebra.

\bigskip

\begin{definition}
A trialgebra $(A,*,\Delta ,\cdot )$ with $*$ and $\cdot $ associative
products on a vector space $A$ (where $*$ may be partially defined, only)
and $\Delta $ a coassociative coproduct on $A$ is given if both $(A,*,\Delta
)$ and $(A,\cdot ,\Delta )$ are bialgebras and the following compatibility
condition between the products is satisfied for arbitrary elements $%
a,b,c,d\in A$: 
\[
(a*b)\cdot (c*d)=(a\cdot c)*(b\cdot d)
\]
whenever both sides are defined.
\end{definition}

\bigskip

Trialgebras were first suggested in \cite{CF} as an algebraic means for the
construction of four dimensional topological quantum field theories. It was
observed there that the representation categories of trialgebras have the
structure of so called Hopf algebra categories (see \cite{CF}) and it was
later shown explicitly in \cite{CKS} that from the data of a Hopf category
one can, indeed, construct a four dimensional topological quantum field
theory. The first explicit examples of trialgebras were constructed in \cite
{GS1} and \cite{GS2} by applying deformation theory, once again, to the
function algebra on the Manin plane and some of the classical examples of
quantum algebras and function algebras on quantum groups. In \cite{GS3} it
was shown that one of the trialgebras constructed in this way appears as a
symmetry of a two dimensional spin system. Besides this, the same trialgebra
can also be found as a symmetry of a certain system of infinitely many
coupled $q$-deformed harmonic oscillators.

The three dimensional extended topological quantum field theories of \cite
{KL} always uniquely determine a finite dimensional trialgebra where one of
the bialgebras contained in it is a quasitriangular Hopf algebra. We call
such trialgebras finite dimensional quasitriangular antipodal trialgbras. On
the other hand, we have the following result:

\bigskip

\begin{lemma}
A finite dimensional quasitriangular antipodal trialgebra always uniquely
determines an algebra object in a modular category.
\end{lemma}

\proof%
Let $\left( T,\cdot ,*,\Delta \right) $ be the trialgebra where $\left(
T,\cdot ,\Delta \right) $ is a quasitriangular Hopf algebra and $\mathcal{C}$
the category of finite dimensional representations of $\left( T,\cdot
,\Delta \right) $. Then $\mathcal{C}$ is modular. Let $A$ be the object in $%
\mathcal{C}$ which corresponds to the defining representation. Then the
product $*$ defines the structure of an algebra object in $\mathcal{C}$ on $%
A $. The product $*$ is an algebra morphism in $\mathcal{C}$ because of the
compatibilities of a trialgebra. 
\endproof%

\bigskip

Observe that the inner algebra objects of modular categories arising from
two dimensional boundary conformal field theories in \cite{FS}, \cite{FRS
2001}, \cite{FRS 2003} do not necessarily determine a trialgebra since the
coproduct leading to the tensor product of $\mathcal{C}$ and the inner
algebra product can, in general, not be realized on one and the same object
of $\mathcal{C}$. So, it seems that the formulation of two dimensional
boundary conformal field theories in terms of modules over inner Frobenius
algebras of modular categories should be a more general framework than one
in terms of trialgebras or of three dimensional extended topological quantum
field theories. Besides this, the approach through modules over inner
Frobenius algebras has the merit that it has been verified rigorously that
nontrivial two dimensional boundary conformal field theories do, indeed,
satisfy its axioms. Since we have found a universal symmetry for trialgebras
and three dimensional extended topological quantum field theories in \cite
{Sch 2002a} and \cite{Sch 2002b}, it is therefore tempting to ask if a
universal symmetry can also be established for pairs $\left( \mathcal{C}%
,M\right) $ consisting of a modular category $\mathcal{C}$ and a module $M$
over an inner Frobenius algebra of $\mathcal{C}$.

Let us now briefly introduce the idea of a universal symmetry of an
algebraic structure. We will restrict to give the general idea, here. For
the needed technical details, we refer the reader to \cite{Dri} (or \cite{CP}
for a short introduction). Consider a quasitriangular quasi-Hopf algebra
with braiding $R$ and Drinfeld associator $\alpha $. In \cite{Dri} Drinfeld
asked for the possibility to change the data $\left( R,\alpha \right) $
while keeping the rest of the structure of the quasitriangular quasi-Hopf
algebra fixed. The Grothendieck-Teichm\"{u}ller group $GT$ is then defined
as the group of such ``gauge transformations'' of the quasitriangular
quasi-Hopf algebra where a certain closure is taken by allowing for formal
deformations of the data $\left( R,\alpha \right) $ in the sense of certain
formal power series. One can show that it is the Drinfeld associator $\alpha 
$ which is the relevant part of the data to basically determine the
structure of $GT$ (see \cite{Dri}, \cite{Kon 1999}). So, we will restrict to
the consideration of the Drinfeld associator, in the sequel.

In \cite{Sch 2002a} we have considered the question of a corresponding
universal symmetry for weak trialgebras where the coproduct receives a
Drinfeld associator $\alpha $ and one of the product a dual Drinfeld
associator $\beta $. We have shown there that the transformations of these
data lead to a self-dual, noncommutative, and noncocommutative Hopf algebra $%
\mathcal{H}_{GT}$ which is a sub-Hopf algebra of the Drinfeld double $%
\mathcal{D}\left( GT\right) $ of $GT$.

\bigskip

Now, let $\mathcal{C}$ be a modular category, $A$ an inner algebra object of 
$\mathcal{C}$, and $M$ a module over $A$, i.e. on $M$ we have a
representation of the algebra product of $A$. We have the folllowing result,
then:

\bigskip

\begin{lemma}
There exists a self-dual, noncommutative, and noncocommutative Hopf algebra $%
\mathcal{H}_F$ which gives the universal symmetry on the pairs $\left( 
\mathcal{C},M\right) $. Besides this, $\mathcal{H}_F$ is a sub-Hopf algebra
of the Drinfeld double $\mathcal{D}\left( GT\right) $ of the
Grothendieck-Teichm\"{u}ller group $GT$ and the Hopf algebra $\mathcal{H}%
_{GT}$ introduced in \cite{Sch 2002a} is a sub-Hopf algebra of $\mathcal{H}%
_F.$
\end{lemma}

\proof%
Following the approach of \cite{Dri} and the introduction of $\mathcal{H}%
_{GT}$ in \cite{Sch 2002a}, we observe that we have two types of data for a
pair $\left( \mathcal{C},M\right) $ the ``gauge freedom'' of which should
determine the universal symmetry: We have the Drinfeld associator $\alpha $
of the modular category $\mathcal{C}$ and we can introduce a similar dual
associator $\beta $ for the representation of the algebra product of $A$ on $%
M$. So, as in \cite{Sch 2002a}, we have, again, a pair $\left( \alpha ,\beta
\right) $ of a Drinfeld associator and a dual Drinfeld associator. Precisely
following the argument given there, we see that the universal symmetry of
the pair $\left( \mathcal{C},M\right) $ has to be given by a subspace 
\[
\mathcal{H}_F\subseteq \mathcal{D}\left( GT\right) 
\]
because if we would transform both, $\alpha $ and $\beta $, separately,
without considering any compatibility constraint between them, each would be
transformed by $GT$. The subspace $\mathcal{H}_F$ takes into account the
constraint that we want to transform a pair $\left( \mathcal{C},M\right) $
into a pair $\left( \mathcal{C},M\right) $.

One proves by calculation that $\mathcal{H}_F$ is closed with respect to
product and coproduct of $\mathcal{D}\left( GT\right) $, i.e. it is a
sub-Hopf algebra of $\mathcal{D}\left( GT\right) $.

Since $\beta $ is a dual Drinfeld associator, i.e. the class of all data $%
\left( \alpha ,\beta \right) $ is self-dual (by $\Bbb{C}$-linearity of $%
\mathcal{C}$ there can not appear more than duals of Drinfeld associators
for $\beta $ and by universality of $\mathcal{H}_F$ all duals of Drinfeld
associators have to appear), and $\mathcal{H}_F$ gives the universal
symmetry, $\mathcal{H}_F$ is self-dual. Using the fact that $GT$ is
non-abelian, we derive that $\mathcal{H}_F$ is noncommutative. Self-duality
then implies that $\mathcal{H}_F$ is also noncocommutative.

Finally, as we mentioned above, the constraint imposed in the definition of
a trialgebra is stronger than the one for an inner algebra object of a
modular category. Therefore $\mathcal{H}_{GT}$ has to be a sub-Hopf algebra
of $\mathcal{H}_F$. 
\endproof%

\bigskip

In \cite{Sch 2002a} we have shown that trialgebras have a strong stability
property which we called ultrarigidity, there. Loosely speaking, trialgebras
can not nontrivially be deformed into algebraic structures with two
associative products and two coassociative coproducts, all linked in a
pairwise compatible way. We have a similar property for pairs $\left( 
\mathcal{C},M\right) $.

\bigskip

\begin{definition}
A bialgebra object in a modular category $\mathcal{C}$ is an object $A$ in $%
\mathcal{C}$ together with two morphisms 
\[
m:A\otimes A\rightarrow A
\]
and 
\[
\Delta :A\rightarrow A\otimes A
\]
in $\mathcal{C}$ where $m$ is associative and $\Delta $ is coassociative and
a morphism of the product $m$.
\end{definition}

\bigskip

\begin{lemma}
Let $\mathcal{C}$ be a bounded modular category and $A$ an algebra object in 
$\mathcal{C}$. Then $A$ can not be nontrivially deformed into a bialgebra
object of a modular category $\widetilde{\mathcal{C}}$.
\end{lemma}

\proof%
Completely similar to the one given in \cite{Sch 2002a}. 
\endproof%

\bigskip

As a consequence, a module $M$ over $A$ can not be nontrivially deformed
into a module over a bialgebra object in some modular category $\widetilde{%
\mathcal{C}}$. Once again, we call this property ultrarigidity.

Observe that nontriviality of the deformation means, here, that the deformed
object is not isomorphic to an algebra object in any modular category. This
does not exclude the possibility of nontrivial deformations with respect to
a fixed modular category, as we know well from the theory of quantum groups.

Since Hopf algebras and modular categories alone do not have such a strong
stability property, this means in the language of physics that two
dimensional boundary conformal field theories are a much more rigid class of
structures than the corresponding bulk theories.

\bigskip

In analogy to the finite dimensional representations of $GT$ (which are
mixed motives of a special kind, see e.g. \cite{Kon 1999}), we will call the
finite dimensional representations of $\mathcal{H}_F$ which are at the same
time corepresentations of $\mathcal{H}_F$ $\mathcal{H}_F$-motives. The
reader, not acquainted with the notions of periods and motives, may think of
a motive as an algebraic abstraction formalizing the period matrices of
algebraic varieties we recall that a period matrix is defined by evaluating
a basis of cohomology on a basis of homology cycles. The coefficients of
such period matrices are called periods. A reader with physics background
may imagine the special periods of representations of $GT$ as the weights of
the Feynman diagrams which appear for the conformal field theories
describing deformation quantization (see \cite{CaFe}, \cite{Kon 1997}).

We have the following result giving a deep interplay between $\mathcal{H}%
_{GT}$ and $\mathcal{H}_F$.

\bigskip

\begin{lemma}
Assume that the conjecture on the transitive action of $GT$ on the algebra $%
P_{\Bbb{Z},Tate}$ of periods of mixed Tate motives over $Spec\left( \Bbb{Z}%
\right) $ of \cite{Kon 1999} is valid. Then the Hopf algebra $\mathcal{H}%
_{GT}$ acts transitively on the class of $\mathcal{H}_F$-motives where we
speak of a transitive action of $\mathcal{H}_{GT}$ if the combined action of 
$\mathcal{H}_{GT}$ on $\alpha $ and its dual on $\beta $ is transitive for
all pairs $\left( \alpha ,\beta \right) $.
\end{lemma}

\proof%
Let $GT_F$ be the subgroup of those elements of $GT$ which are consistent
with the constraint leading from $\mathcal{D}\left( GT\right) $ to $\mathcal{%
H}_F$, i.e. which can appear as a transformation of one of the components $%
\left( \alpha ,\beta \right) $ such that the transformation is consistent
with a complete transformation of the pair $\left( \alpha ,\beta \right) $.
Define in a similar way $GT_c$ as the subgroup consisting of those elements
of $GT$ which are consistent with the constraint of $\mathcal{H}_{GT}$.
Using self-duality of $\mathcal{H}_{GT}$ and $\mathcal{H}_F$ and 
\[
\mathcal{H}_{GT}\subseteq \mathcal{H}_F 
\]
it follows that 
\[
GT_c\subseteq GT_F 
\]
Making use of the explicit examples of trialgebras constructed in \cite{GS2}%
, one can prove by direct calculation that 
\[
GT_c=GT 
\]
because the tensor product construction used in \cite{GS2} is not affected
by the introduction of a Drinfeld associator for the Hopf algebra one starts
from. So, it follows that 
\[
GT_F=GT 
\]
But then a $\mathcal{H}_F$-motive is always given by a pair of finite
dimensional representations of $GT$ plus a constraint deriving from the
defining constraint of $\mathcal{H}_F$.

But the above mentioned conjecture of \cite{Kon 1999} and a theorem of Nori
which states that the algebra of periods is isomorphic to the algebra of
functions on the torsor of isomorphisms between algebraic de Rham and Betti
cohomology, $GT$ acts transitively on its finite dimensional
representations. But since the constraint for $\mathcal{H}_F$-motives is of
a type which requires that a pair of $GT$-representations defines a special
Hopf algebra, it follows from the definition of $\mathcal{H}_{GT}$ in \cite
{Sch 2002a} that $\mathcal{H}_{GT}$ acts transitively on the $\mathcal{H}_F$%
-motives. 
\endproof%

\bigskip

\begin{remark}
This result shows that though 
\[
\mathcal{H}_{GT}\subseteq \mathcal{H}_F
\]
$\mathcal{H}_{GT}$ acts on the representation theoretic side as a symmetry
for $\mathcal{H}_F$. In a subsequent publication, we plan to put the results
given here into a broader algebraic context. There is some hope that in that
context one can establish the equivalence of $\mathcal{H}_{GT}$ and $%
\mathcal{H}_F$ if using a suitable notion of Morita equivalence for
self-dual Hopf algebras.
\end{remark}

\bigskip

In the next section, we will make use of the above result in order to derive
a transitive action of a sub-Hopf algebra $\mathcal{H}_D$ of $\mathcal{H}%
_{GT}$ on the relative period matrices of holomorphic $N=1$ special
geometries.

\bigskip

\section{Application in holomorphic $N=1$ special geometry}

For a type II string theory compactified on a three complex dimensional
Calabi-Yau manifold $Y$ one can show that the corresponding string field
theory in the so called $B$-model is given by the Kodaira-Spencer equation
of holomorphic deformations of a complex manifold (see \cite{BCOV}). One can
further show that holomorphic deformations in\ Kodaira-Spencer theory act on
the third cohomology $H^3\left( Y\right) $ of $Y$ and this action completely
specifies the deformation theory. Concretely, the period matrix of $Y$
completely determines the structure constants of the operator product
expansion of the $B$-model with target space $Y$ (see \cite{BCOV} and for
short overviews the introductions of \cite{LMW 2002a}, \cite{LMW 2002b}).

The papers \cite{LMW 2002a} and \cite{LMW 2002b} deal with the question if
this approach can be extended to the case of two dimensional boundary
conformal field theory, specifically to a $B$-model of open-closed strings.
The answer is that one gets relative Kodaira-Spencer theory, i.e. the
holomorphic deformation theory of $Y$ together with a complex submanifold $%
X\subseteq Y$ such that both, $X$ and $Y$, receive holomorphic deformations
but the submanifold property is kept. The answer found there is that
relative Kodaira-Spencer theory leads to and is determined by an action on
the third relative cohomology $H^3\left( X,Y\right) $ (see the cited papers
or \cite{KaLe} for an introduction to relative cohomology) and the relative
period matrix allows one to completely determine the structure constants of
the extended chiral ring (which includes the boundary operators).

First, observe that the usual (absolute) period matrices correspond to mixed
motives, i.e. finite dimensional representations of the motivic Galois
group. Using the just mentioned results on relative period matrices and
relative Kodaira-Spencer theory, it follows that a relative period matrix
has to correspond to a pair of mixed motives satisfying a certain constraint
deriving from the submanifold property 
\[
X\subseteq Y
\]
On the other hand, using the result that relative Kodaira-Spencer theory is
precisely the string field theory of the open-closed $B$-model, we know that
this geometric submanifold constraint has to be equivalent to the algebraic
constraint of the previous section which appears from the inner Frobenius
algebra condition in modular categories. It follows that relative period
matrices are restricted to correspond to pairs of matrices determined by the
quotient group $GT$ of the motivic Galois group and that the constraint
precisely restricts them to the type of a $\mathcal{H}_F$-motive.

Using the results of the previous section, we immediately have the following
corollary, then:

\bigskip

\begin{corollary}
There exists a sub-Hopf algebra $\mathcal{H}_D$ of $\mathcal{H}_{GT}$ that
acts transitively on the relative period matrices of holomorphic $N=1$
special geometries (i.e. on the relative period matrices of the open-closed $%
B$-model).
\end{corollary}

\proof%
By the results of the previous section, $\mathcal{H}_{GT}$ acts transitively
on $\mathcal{H}_F$-motives. It follows that there is a sub-Hopf algebra $%
\mathcal{H}_D$ of $\mathcal{H}_{GT}$ which respects the constraint that a $%
\mathcal{H}_F$-motive corresponds to a relative period matrix of holomorphic 
$N=1$ special geometry. Obviously, $\mathcal{H}_D$ acts transitively on the
relative period matrices. 
\endproof%

\bigskip

\begin{remark}
In a later publication, we hope to show that $\mathcal{H}_D$ and $\mathcal{H}%
_{GT}$, again, should be Morita equivalent.
\end{remark}

\bigskip

The above corollary means that $\mathcal{H}_D$ acts on holomorphic $N=1$
special geometries as a duality Hopf algebra in the sense that any two such
geometries are dual to each other with respect to the $\mathcal{H}_D$-action.

\bigskip

\section{Conclusion}

We have shown the existence of a universal symmetry for Frobenius algebras
in modular categories given by a self-dual, noncommutative, and
noncocommutative Hopf algebra $\mathcal{H}_F$. We have discussed the
relationship between $\mathcal{H}_F$ and the Hopf algebra $\mathcal{H}_{GT}$
introduced previously. Finally, we have applied our result to arrive at a
sub-Hopf algebra $\mathcal{H}_D$ of $\mathcal{H}_{GT}$ acting as a duality
Hopf algebra on holomorphic $N=1$ special geometries.

It remains a task for future work to discuss the concrete physical
implications of the appearance of the Hopf algebra $\mathcal{H}_D$ in
holomorphic $N=1$ special geometry. As we have remarked already above, we
plan to put the results gained, here, in a broader algebraic context in a
future publication, in order to achieve a deeper understanding of the
appearance of these Hopf algebraic structures on the moduli spaces of
compactifications in string theory.

\bigskip

\textbf{Acknowledgements:}

I thank the Erwin Schr\"{o}dinger Institute for Mathematical Physics,
Vienna, and the Institute for Theoretical Physics of the University of
Innsbruck for hospitality during part of the time in which this work has
been done. For discussions on the subjects involved, I thank Bojko Bakalov,
Alberto Cattaneo, Karen Elsner, Harald Grosse, Christoph Schweigert, and
Ivan Todorov.


\bigskip

\begin{thebibliography}{LMW 2002a}
\bibitem[BCOV]{BCOV}  M. Bershadsky, S. Cecotti, H. Ooguri, C. Vafa, \textit{%
Kodaira-Spencer theory of gravity and exact results for quantum string
amplitudes}, Comm. Math. Phys. 165 (1994), 311-427.

\bibitem[CaFe]{CaFe}  A. Cattaneo, G. Felder, \textit{A path integral
approach to the Kontsevich quantization formula}, math/9902090.

\bibitem[CF]{CF}  L. Crane, I. B. Frenkel, \textit{Four-dimensional
topological quantum field theory, Hopf categories and the canonical bases},
J. Math. Phys. 35, 5136-5154 (1994).

\bibitem[CKS]{CKS}  J. S. Carter, L. H. Kauffman, M. Saito, \textit{%
Structures and diagrammatics of four dimensional topological lattice field
theories}, math.GT/9806023.

\bibitem[CP]{CP}  V. Chari, A. Pressley, \textit{Quantum groups}, Cambridge
University Press, Cambridge 1994.

\bibitem[Dri]{Dri}  V. G. Drinfeld, \textit{On quasi-triangular Quasi-Hopf
algebras and a group closely related with Gal(}$\overline{\Bbb{Q}}/\Bbb{Q}$%
\textit{)}, Leningrad Math. J., \textbf{2}, 829-860 (1991).

\bibitem[FS]{FS}  J. Fuchs, C. Schweigert, \textit{Category theory for
conformal boundary conditions}, math.CT/0106050v2.

\bibitem[FRS 2001]{FRS 2001}  J. Fuchs, I. Runkel, C. Schweigert, \textit{%
Conformal correlation functions, Frobenius algebras and triangulations},
hep-th/0110133v2.

\bibitem[FRS 2003]{FRS 2003}  J. Fuchs, I. Runkel, C. Schweigert, \textit{%
Boundaries, defects and Frobenius algebras}, hep-th/0302200.

\bibitem[GS 2000a]{GS1}  H. Grosse, K.-G. Schlesinger, \textit{On a
trialgebraic deformation of the Manin plane}, Letters in Mathematical
Physics 52 (2000), 263-275.

\bibitem[GS 2000b]{GS2}  H. Grosse, K.-G. Schlesinger, \textit{On second
quantization of quantum groups}, J. Math. Phys. 41 (2000), 7043-7060,
earlier version also available in the preprint series of the Erwin
Schr\"{o}dinger Institute for Mathematical Physics, Vienna, 841 (under
http://www.esi.ac.at).

\bibitem[GS2001]{GS3}  H. Grosse, K.-G. Schlesinger, \textit{A suggestion
for an integrability notion for two dimensional spin systems}, Letters in
Mathematical Physics 55 (2001), 161-167, hep-th/0103176.

\bibitem[KaLe]{KaLe}  M. Karoubi, C. Leruste, \textit{Algebraic topology via
differential geometry}, Cambridge University Press, Cambridge 1987.

\bibitem[KL]{KL}  T. Kerler, V. V. Lyubashenko, \textit{Non-semisimple
topological quantum field theories for 3-manifolds with corners}, Lecture
Notes in Mathematics 1765, Springer, Berlin and Heidelberg 2001.

\bibitem[Kon 1997]{Kon 1997}  M. Kontsevich, \textit{Deformation
quantization of Poisson manifolds I}, math/9709180.

\bibitem[Kon 1999]{Kon 1999}  M. Kontsevich, \textit{Operads and motives in
deformation quantization}, math.QA/9904055.

\bibitem[LMW 2002a]{LMW 2002a}  W. Lerche, P. Mayr, N. Warner, \textit{%
Holomorphic N=1 special geometry of open-closed type II strings},
hep-th/0207259.

\bibitem[LMW 2002b]{LMW 2002b}  W. Lerche, P. Mayr, N. Warner, \textit{N=1
special geometry, mixed Hodge variations and toric geometry}, hep-th/0208039.

\bibitem[Sch 2002a]{Sch 2002a}  K.-G. Schlesinger, \textit{A quantum
analogue of the Grothendieck-Teichm\"{u}ller group}, J. Phys. A: Math. Gen.
35 (2002), 10189-10196, math.QA/0104275.

\bibitem[Sch 2002b]{Sch 2002b}  K.-G. Schlesinger, \textit{A universal
symmetry structure in open string theory}, hep-th/0203183.
\end{thebibliography}
\end{document}